# Vectored route-length minimization – a heuristic and an open conjecture

Florentin Smarandache

The University of New Mexico, USA

Sukanto Bhattacharya

The University of Queensland, Australia

**Abstract**

We have posed a simple but interesting graph theoretic problem and posited a heuristic solution procedure, which we have christened as <u>*Ve*</u>ctored <u>*R*</u>oute-length <u>*M*</u>inimization <u>*S*</u>earch (VeRMinS). Basically, it constitutes of a re-casting of the classical "shortest route" problem within a strictly Euclidean space. We have only presented a heuristic solution process with the hope that a formal proof will eventually emerge as the problem receives wider exposure within mathematical circles.

**A short historical background of similar class of problems**

The classical "shortest route" (or shortest path) problem is properly associated with the branch of mathematics formally known as *graph theory* (or *network theory*). History has it that this theory originated in an attempt to solve a famous 18$^{th}$ century routing problem concerning the Prussian city of Konigsberg (Kaliningrad in modern Russia). The city is located along the two banks and on two islands formed by the river Pregel, which effectively divides the city into four separate landmasses. Seven bridges connected the various regions of the city and the resulting "Konigsberg bridge problem" had to do with finding an optimal route around the city that would require a traveler to cross each of the seven bridges only once in the whole trip (Alexanderson, 2006). That it is impossible to make such a trip was originally proved by the Swiss mathematical genius Leonhard Euler (Euler, 1766) thus formally giving birth to the mathematics of networks.

The classical "shortest route" problem born out of the Konigsberg bridge problem subsequently branched into a number of well-known variants popularly grouped as "traveling salesman" problems. The shortest route problem is one of the many practical adaptations of Eulerian graph theory. The basic problem is concerned with finding the shortest distance between a "source" and a "sink" node in any sufficiently generalized network consisting of a finite number of nodes.

The usual practical applications of similar class of problems in modern times are in the configuration of telecommunications networks e.g. connecting one transmission tower to another in a network so that the total network up-linking time is minimized. There are also interesting application possibilities in the realm of social sciences especially in *social network analysis* that has provided valuable insights into the governance and behavior of organized groups in society and social capital generation (Nan Lin, 1999). In business and finance applications, network data mining is being applied to detect fraud and money laundering activities (Yue et. al., 2007) and in following terrorist money trails by identifying the likely "shortest paths" through social networks (Keefe, 2006)

Many alternative algorithms to solving the shortest route problem have been devised e.g. *Djikstra's algorithm* (Djikstra, 1959), and *Ford-Fulkerson's algorithm* (Ford and Fulkerson, 1962), which have many applications in the fields of telecommunications and internetworking. However, in positing our problem, we have been concerned with the most simplistic version of the classical shortest route problem in strictly Euclidean space of unrestricted dimensionality, which we proceed to define as follows:

"*Given a partially connected network of N nodes in a strictly Euclidean space of any dimension, find a route through the network from a pre-specified source node $S_0$ to a pre-specified sink node $S_N$ such that the overall route length (in terms of the total Euclidean distance) is minimum*"

**Mathematical basis of VeRMinS: a proposed heuristic solution procedure**

*The Vectored Route-length Minimization Search* (VeRMinS) is a *heuristic search* that aims to find the shortest route from a source node to a sink node in a network in Euclidean space of any dimension by identifying the *linear-most connectivity* between the source and sink nodes.

With every route in a network, we associate a corresponding *weight* factor, which is the sum of the Euclidean distance between the nodes on that route. Then the best (i.e. linear-most) route through the network is the one having the minimum weight (Rote, 1990). For any network consisting of $N = m + 1$ nodes, we can set up a *network connectivity matrix* **M** as follows:

| 1 | $R_{01}$ | $R_{02}$ | … | $R_{0m}$ |
|---|---|---|---|---|
| $R_{10}$ | 1 | $R_{12}$ | … | $R_{1m}$ |
| $R_{20}$ | $R_{21}$ | 1 | … | $R_{2m}$ |
| … | … | … | … | … |
| $R_{m0}$ | $R_{m1}$ | $R_{m2}$ | … | 1 |

In the network connectivity matrix, when $i \neq j$, $R_{ij} = 1$ if and only if a connectivity exists between nodes i and j and $R_{ij} = 0$ otherwise. Since a node is necessarily 'self-connected', $R_{ij} = 1$ when $i = j$ i.e. for all the diagonal elements of **M**.

A finite number, say q, of *route vectors* $\mathbf{P}_t$ (with $t = 1, 2 \ldots q$) can then be extricated from **M** such that $\mathbf{P}_1 = [k_{10}\ k_{11}\ k_{12}\ \ldots\ k_{1m}]$, $\mathbf{P}_2 = [k_{20}\ k_{21}\ k_{22}\ \ldots\ k_{2m}]\ \ldots\ \mathbf{P}_t = [k_{t0}\ k_{t1}\ k_{t2}\ \ldots\ k_{tm}]\ \ldots\ \mathbf{P}_q = [k_{q0}\ k_{q1}\ k_{q2}\ \ldots\ k_{qm}]$, where $k_{qj} = 1$ if node j lies on the q-th route and $k_{qj} = 0$ otherwise.

A (m x 1) weight vector **W** is defined as follows:

$$\mathbf{W} = [0\ w_1\ w_2\ \ldots\ w_i,\ w_{m-1}\ 0]^T$$

where $w_i$ is the vertical Euclidean distance of the i-th node from the *ideal route* (which is simply a hypothetical straight line connecting the source and sink nodes), as determined by its position vector with respect to the *ideal route*. Since both the source and sink nodes must necessarily lie on the shortest route (i.e. a route must be *effective* before it can be *efficient*), $w_0 = w_m = 0$.

Then, $\mathbf{P}_q \cdot \mathbf{W} = \sum\sum (k_{qj} w_i)$ would yield the deciding criterion for the q-th route in terms of the vertical Euclidean distances of each of the nodes along the q-th route from the *ideal* route.

**Introducing the property of *Euclidean dominance***

*The route vector $\mathbf{P}_a$ exhibits Euclidean dominance over the route vector $\mathbf{P}_b$ (written henceforth as $\mathbf{P}_a \supset \mathbf{P}_b$) when at least one element of $\mathbf{P}_a$ is 0 with the corresponding element in $\mathbf{P}_b$ being 1 and all other elements being same for $\mathbf{P}_a$ and $\mathbf{P}_b$.*

**Proof:** This property follows from the principle of triangular inequality in Euclidean geometry whereby the sum of two sides of a triangle is always greater in magnitude than the third side.

Each of the nodes in a network corresponds to a particular position vector in Euclidean space. Therefore, it implies that if node A is connected to both nodes B and C while node B is also connected to node C, then the route that goes directly from node A to node C will always be more preferable than one which goes from node A to node B to node C. This of course assumes that the remaining segments of the two routes coincide with each other.

So the property of Euclidean dominance may be used to effectively eliminate some of the q route vectors extricated from $\mathbf{M}$. Assuming h route vectors are eliminated after applying Euclidean dominance, then the *linear-most route* is obtainable as $\text{Min}_t\ [\mathbf{P}_1 \cdot \mathbf{W},\ \mathbf{P}_2 \cdot \mathbf{W},\ \ldots,\ \mathbf{P}_t \cdot \mathbf{W},\ \ldots,\ \mathbf{P}_{(q-h)} \cdot \mathbf{W}]$.

**Applying the VeRMinS – a numerical illustration**

Let a simple network in 2*D*-Euclidean space consisting of ten nodes 0, 1, 2 … 9 be as follows:

| Preceding node | Succeeding node | $w_i$ |
|---|---|---|
| 0 | 1, 2, 3 | 0 |
| 1 | 4, 7 | 3 |
| 2 | 4, 5, 6 | 0 |
| 3 | 6, 8 | 3 |
| 4 | 7 | 2 |

| 5 | 7, 8 | 0 |
| 6 | 8 | 1 |
| 7 | 9 | 5 |
| 8 | 9 | 6 |
| 9 | - | 0 |

We wish to find the best (i.e. linear-most) route from node 0 to node 9.

We wish to make the readers aware that here we only present an illustrative exercise outlining a numerical solution procedure. However, we supply no formal proof that the outlined procedure is necessary and sufficient in obtaining the shortest route through a network in any Euclidean space.

The network connectivity matrix $\mathbf{M}_{10 \times 10}$ for the above network is obtained as follows:

| 1 | 1 | 1 | 1 | 0 | 0 | 0 | 0 | 0 | 0 |
|---|---|---|---|---|---|---|---|---|---|
| 1 | 1 | 0 | 0 | 1 | 0 | 0 | 1 | 0 | 0 |
| 1 | 0 | 1 | 0 | 1 | 1 | 1 | 0 | 0 | 0 |
| 1 | 0 | 0 | 1 | 0 | 0 | 1 | 0 | 1 | 0 |
| 0 | 1 | 0 | 0 | 1 | 0 | 0 | 1 | 0 | 0 |
| 0 | 0 | 1 | 0 | 0 | 1 | 0 | 1 | 1 | 0 |
| 0 | 0 | 1 | 1 | 0 | 0 | 1 | 0 | 1 | 0 |
| 0 | 1 | 0 | 0 | 1 | 1 | 0 | 1 | 0 | 1 |
| 0 | 0 | 0 | 1 | 0 | 1 | 1 | 0 | 1 | 0 |

| 0 | 0 | 0 | 0 | 0 | 0 | 0 | 1 | 1 | 1 |

The following route vectors may be extricated from **M**:

**P**$_1$ = [1 1 0 0 1 0 0 1 0 1]

**P**$_2$ = [1 1 0 0 0 0 0 1 0 1]

**P**$_3$ = [1 0 1 0 1 0 0 1 0 1]

**P**$_4$ = [1 0 1 0 0 1 0 1 0 1]

**P**$_5$ = [1 0 1 0 0 1 0 0 1 1]

**P**$_6$ = [1 0 1 0 0 0 1 0 1 1]

**P**$_7$ = [1 0 0 1 0 0 1 0 1 1], and

**P**$_8$ = [1 0 0 1 0 0 0 1 1].

It may be easily observed that **P**$_2$ ⊃ **P**$_1$ and **P**$_8$ ⊃ **P**$_7$, so, using the property of Euclidean dominance one can eliminate **P**$_1$ and **P**$_7$ straightaway.

The weight vector is obtained as: **W** = [0 3 0 3 2 0 1 5 6 0]$^T$

Therefore **P**$_2$.**W** = 8, **P**$_3$.**W** = 7, **P**$_4$.**W** = 5, **P**$_5$.**W** = 6, **P**$_6$.**W** = 7 and **P**$_8$.**W** = 9.

So **W**\* = Min$_t$ [**P**$_t$.**W**] = 5, which corresponds to the route vector **P**$_4$ thereby identifying it as the linear-most route from source to sink.

**An open conjecture**

VeRMinS is proposed at this stage as no more than a heuristic search procedure. We have not supplied a formal proof that the outlined search procedure is necessary and sufficient in obtaining the shortest route through a network of a finite number of nodes in any Euclidean space of unrestricted dimensionality. This problem is left open at this stage that may either be proved by showing that all other possible search procedures will always yield less optimal (i.e.

longer) routes or disproved via a counter-example that shows that a shorter route exists through a network in any strictly Euclidean space that is not picked by the outlined VeRMinS procedure.